\documentclass{amsart}
\usepackage{amssymb,mathrsfs,amsfonts}
\newtheorem{pro}{Proposition}[section]
\newtheorem{nota}{Remark}[section]
\newtheorem{note}{Remarks}[section]
\title{On Generalized Cauchy-Stieltjes transforms of some Beta distributions}

\begin{document}
\maketitle
\centerline{Nizar Demni\footnote{Fakult\"at f\"ur mathematik, SFB  701, universit\"at Bielefeld, Bielefeld, Germany, email: demni@math.uni-bielefeld.de.}}

\begin{abstract}
We express the generalized Cauchy-Stieltjes transforms (GCST) of some particular Beta distributions depending on a positive parameter $\lambda$ as $\lambda$-powered Cauchy-Stieltjes transforms (CST) of some probability measures. The CST of the latter measures are shown to be the geometric mean of the CST of the Wigner law together with another one. Moreover, they are absolutely continuous and we derive their densities by proving that they are the so-called Markov transforms of compactly-supported probability distributions. Finally, a detailed analysis is performed on one of the symmetric Markov transforms which interpolates between the Wigner ($\lambda = \infty$) and the arcsine
($\lambda = 1$) distributions. We first write down its moments through a terminating series ${}_3F_2$ and show that they are polynomials in the variable $1/\lambda$, however they are no longer positive integer-valued as for $\lambda = 1, = \infty$ (for instance $\lambda=2$) thereby no general combinatorial interpretation holds. Second, we compute the free cumulants in the case when $\lambda = 2$ and explain how to proceed in the cases when $\lambda = 3,4$. Problems of finding a parallel to the representation theory of the infinite symmetric group and an interpolating convolution are discussed.    
\end{abstract}

\section{Motivation}
Let $\lambda > 0$ and $\mu_{\lambda}$ a probability measure (possibly depending on $\lambda$) with finite all order moments. The \emph{generalized Cauchy-Stieltjes transform} (GCST) of $\mu_{\lambda}$ is defined by
\begin{equation*}
\int_{\mathbb{R}} \frac{1}{(z-x)^{\lambda}}\mu_{\lambda}(dx) 
\end{equation*}
for non real complex $z$ (\cite{Sum}, \cite{Hirs}, \cite{Sch}). For $\lambda = 1$, it reduces to the (ordinary) \emph{Cauchy-Stieltjes transform} (CST) which has been of great importance during the last two decades for both probabilists and algebraists due the central role it plays in free probability and representation theories  (\cite{Bia}, \cite{Ker}). Moreover, CST were extensively studied and they are well-understood in the sense that for instance, a complete characterization of those functions is known and one has a relatively easy inversion formula due to Stieltjes (\cite{Ker}). However, their generalized versions are more hard to handle as one may realize from the complicated inversion formulas in \cite{Sum}, \cite{Hirs}, \cite{Sch}. 
In this paper, we adress the problem of relating both generalized and ordinary transforms, that is, given $\mu_{\lambda}$ check whether there exists a {\bf probability} measure $\nu_{\lambda}$ such that 
\begin{equation}\label{M0}
\int_{\mathbb{R}} \frac{1}{(z-x)^{\lambda}}\mu_{\lambda}(dx) = \left[\int_{\mathbb{R}} \frac{1}{z-x}\nu_{\lambda}(dx)\right]^{\lambda}
\end{equation}  
for $z$ in some suitable complex region where $(z-x)^{\lambda}$ is single-valued and is analytic (but we prefer taking the principal determination) and characterize 
$\nu_{\lambda}$ in the affirmative case. Doing so will lead for instance to the invertibility of 
\begin{equation*}
z \mapsto \left[\int_{\mathbb{R}} \frac{1}{(z-x)^{\lambda}}\mu_{\lambda}(dx)\right]^{1/\lambda} 
\end{equation*}
for $z$ belonging to some neighborhood of infinity (\cite{Ber}) and as a by-product to a kind of $\lambda$-free cumulants generating function for $\mu_{\lambda}$, referring to the case $\lambda =1$ (\cite{Spei}). 
However, this may not be always possible as we shall see, that is, $\nu_{\lambda}$ may not be a probability measure for some $\lambda$. Moreover, when $\nu_{\lambda}$ is shown to be a probability measure, it is not easy to check whether it is absolutely continuous or not and more harder will be to write down its density when it is so. This was behind our willing to get an insight into the above problem since GCST naturally appear in the study of a particular class of probability measures $\mu_{\lambda}$ 
(\cite{Dem}). Those probability measures may be mapped via affine transformations into Beta distributions 
\begin{equation*}
\beta_{a,b}(dx) = (2-x)^{a}(2+x)^{b} {\bf 1}_{[-2,2]}(x) dx
\end{equation*} 
with parameters $a,b > -1$ depending on $\lambda$, and their monic orthogonal polynomials, say $(P_n^{\lambda})_n$, are the only ones that admit an ultraspherical type 
generating function: 
\begin{equation}\label{GF}
\sum_{n \geq 0}\frac{(\lambda)_n}{n!} P_n^{\lambda}(x) z^n = \frac{1}{u_{\lambda}(z)(f_{\lambda}(z) -x)^{\lambda}}, \quad x \in \textrm{supp}(\mu_{\lambda}),
\end{equation}        
where $u_{\lambda}, f_{\lambda}$ satisfy some technical conditions in a suitable open complex region near $z=0$ so that (\ref{GF}) makes sense.
Due to the orthognality of $P_n^{\lambda}, n \geq 0$, one gets after integrating both sides in (\ref{GF}) with respect to $\mu_{\lambda}$ that 
 \begin{equation*}
u_{\lambda}(z) = \int_{\mathbb{R}}\frac{1}{(f_{\lambda}(z)- x)^{\lambda}}\mu_{\lambda}(dx),
\end{equation*}
for small $z$. Thus, if $f_{\lambda}$ is invertible, the last equality is rewritten as 
\begin{equation*}
u_{\lambda}[f_{\lambda}^{-1}(z)] = \int_{\mathbb{R}}\frac{1}{(z- x)^{\lambda}}\mu_{\lambda}(dx)
\end{equation*}
so that $u_{\lambda}(f_{\lambda}^{-1})$ is the GCST of $\mu_{\lambda}$. Fortunately, as one easily sees from the reminder below, $f_{\lambda}^{-1}$ may be identified with the 
CST (in some neighborhood of infinity) of a semi-circle law of mean and variance depending on $\lambda$. 

\section{Reminder and results}
Assume $\mu_{\lambda}$ has zero mean and unit variance, then $\mu_{\lambda}$ has a generating function for orthogonal polynomials of ultraspherical-type if and only if it is compactly-supported and belongs to one of the four families given by the following data (\cite{Dem}): 
\begin{equation*}
\left.\begin{array}{ll}
u_{\lambda}(z) = z^{\lambda}, & f_{\lambda}(z) = \displaystyle \frac{\lambda+1}{2}z + \frac{1}{z},  \lambda > 0, \\

u_{\lambda}(z) = \displaystyle \frac{z^{\lambda}}{1- (\lambda/2)z^2}, & f_{\lambda}(z) = \displaystyle \frac{\lambda}{2}z + \frac{1}{z}, \lambda > 1/2 \\

u_{\lambda}(z) = \displaystyle \frac{z^{\lambda}}{1\pm \lambda z/\sqrt{2\lambda-1}}, & f_{\lambda}(z) = \displaystyle \frac{\lambda^2}{2\lambda-1}z \pm \frac{1}{\sqrt{2\lambda-1}} + \frac{1}{z}, \lambda > 1/2. 
\end{array}\right.
\end{equation*}                     
We did not displayed the density of $\mu_{\lambda}$ for sake of clarity and we will do it later. Nevertheless, it is worth recalling that the first family corresponds to the monic Gegenbauer polynomials $(C_n^{\lambda})_n$, the second one corresponds to $(C_n^{\lambda-1})_n$ and is related to the Poisson kernel, while the remaining families correspond to shifted Jacobi polynomials whose paremeters differ by $1$. Now, the inverse in composition's sense of the CST $G_{m,\sigma}$ of a semi-circular law of mean $m$ and variance $\sigma$, say $K_{m,\sigma}$, is given by (\cite{Boz})
\begin{equation}\label{K-Tr}
K_{m,\sigma}(z) := G^{-1}_{m,\sigma}(z) = \sigma^2z + a+ \frac{1}{z}. 
\end{equation} 
Then, the reader can check our previous claim about $f_{\lambda}$. Henceforth, we will only make use of $G:= G_{0,1}$ since the elemantary identity holds 
\begin{equation*}
G_{m,\sigma}(z) = \frac{1}{\sigma}G\left(\frac{x-m}{\sigma}\right).
\end{equation*}
The first main result may be stated for sufficiently large $z$ as 
\begin{equation}\label{E1}
\left[\int_{\mathbb{R}}\frac{1}{(z- x)^{\lambda}}\beta_{a,b}(dx)\right]^{1/\lambda} = [\tilde{u}_{\lambda}(G(z))]^{1/\lambda} = G^{\alpha(\lambda)}(z)\tilde{G}(z)^{\gamma(\lambda)}, 
\end{equation} 
where $\alpha(\lambda)  + \gamma(\lambda) = 1$ and
\begin{equation*}
\left.\begin{array}{ll}
\tilde{u}_{\lambda}(z) = z^{\lambda}, & \alpha(\lambda) = 1,   \\

\tilde{u}_{\lambda}(z) = \displaystyle \frac{z^{\lambda}}{1- z^2}, & \alpha(\lambda) = 1- \displaystyle \frac{1}{\lambda}, \tilde{G}(z) = \frac{1}{\sqrt{z^2-4}}, \\

\tilde{u}_{\lambda}(z) = \displaystyle \frac{z^{\lambda}}{1\pm z}, &  \alpha(\lambda) = 1- \displaystyle \frac{1}{2\lambda}, \tilde{G}(z) = \frac{1}{z\pm2}
\end{array}\right.
\end{equation*}                     
for the four families respectively. It follows by the characterization of CST of probability measures that 
\begin{equation*}
z\ \mapsto \left[\int_{\mathbb{R}}\frac{1}{(z- x)^{\lambda}}\beta_{a,b}(dx)\right]^{1/\lambda}  
\end{equation*} 
defines a CST of some probability measure $\nu_{\lambda}$ for $\lambda$ provided that $0 \leq \alpha(\lambda) = 1 - \gamma(\lambda) \leq 1$. Note that under this condition, the CST of $\nu_{\lambda}$ is the geometric mean of $G,\tilde{G}$ and that one discards the values $\lambda \in ]1/2,1[$ for the second family for which $\tilde{G}$ is the CST of the arcsine distribution (\cite{Spei}).\\ 
The second main result states that under the same condition, $\nu_{\lambda}$ is an absolutely continuous probability measure and gives its density. This follows from the fact that $\nu_{\lambda}$ is the so-called Markov transform of some compactly-supported probability measure $\tau_{\lambda}$ (\cite{Ker}), that is  
\begin{equation}\label{MI}
\int_{\mathbb{R}} \frac{1}{z-x}\nu_{\lambda}(dx) = \exp -\int_{\mathbb{R}} \log (z-x) \tau_{\lambda}(dx),
\end{equation}
and from Cifarelli and Regazzini's results (\cite{Ker} p. 51). For the first and the second families, $\nu_{\lambda}$ is given by the Wigner distribution and a symmetric deformation of it respectively, while for the remaining ones, it is a non symmetric deformation of the Wigner distribution. Since the latter is a universal limiting object (representation theory of the infinite symmetric group, spectral theory of large random matrices, free probability theory), we give a particular interest in $\nu_{\lambda}$ corresponding to the second family. 
We first express its moments by means of a terminating ${}_3F_2$ series interpolating between the moments of the Wigner and the arcsine distributions (Catalan and shifted Catalan numbers, \cite{Ker} p.64). Unfortunately, the moments are no longer positive integer-valued  as it is shown for $\lambda=2$ therefore no general combinatorial interpretation holds. Nevertheless, they are polynomials in the variable $1/\lambda$. Finally, we use (\ref{E1}) to compute the inverse of its CST and the free cumulants generating function in the case 
$\lambda=2$  which involve a weighted sum of the Catalan and the shifted Catalan numbers. For $\lambda = 3,4$, this is a more complicated task and by Galois theory, this is a limitation rather than a restriction since we are led to find a root of a polynomial equation of degree $\lambda$. 

\begin{note}
1/For the second family and for the discarded values of $\lambda$ such that the condition $0 \leq \alpha(\lambda) = 1 - \gamma(\lambda) \leq 1$ 

\begin{equation*}
z \mapsto \left[\int_{\mathbb{R}} \frac{1}{(z-x)^{\lambda}}\mu_{\lambda}(dx)\right]^{1/\lambda}
\end{equation*}  
does not define the CST of a probability distribution.\\ 
2/ It follows from (\ref{M0}) and (\ref{MI}) that 
\begin{align*}
\int_{\mathbb{R}} \frac{1}{(z-x)^{\lambda}}\mu_{\lambda}(dx) & = \exp - \lambda \int_{\mathbb{R}} \log (z-x) \tau_{\lambda}(dx)
= \left[\int_{\mathbb{R}} \frac{1}{z-x}\nu_{\lambda}(dx)\right]^{\lambda}.
\end{align*}
Similar identities already showed up in relation to Bayesian statistics (\cite{Ker} p. 59).
\end{note}

Throughout the paper, computations are performed up to constants depending on $\lambda$, which normalize the finite positive measures involved here to be probability measures. The paper is divided into five sections: the first four sections are devoted to the four families $\mu_{\lambda}$ while the last one is devoted to the particular interest we give in the probability measure 
$\nu_{\lambda}$ corresponding to the second family. 

\section{Markov transforms: symmetric measures}
\subsection{First family}
On the one hand, 
\begin{equation*}
\mu_{\lambda}(dx)  \propto \left(1- \frac{x^2}{2(1+\lambda)}\right)^{\lambda-1/2} {\bf 1}_{[\pm \sqrt{2(1+\lambda)}]}(x)dx, \lambda > 0,
\end{equation*} 
and its image of under the map $x \mapsto \sqrt{(1+\lambda)/2}x$ has the density proportionaI to 
\begin{equation*}
\left(4- x^2\right)^{\lambda-1/2}{\bf 1}_{[-2,2]}(x).
\end{equation*}
On the other hand, it is easy to see that 
\begin{equation*}
f_{\lambda}(z) = \frac{\lambda+1}{2}z + \frac{1}{z} = K_{0,\sqrt{(1+\lambda)/2}}(z) 
\end{equation*}
so that 
\begin{equation*}
f_{\lambda}^{-1}(z) = G_{0,\sqrt{(1+\lambda)/2}}(z) = \sqrt{\frac{2}{1+\lambda}} G\left(\sqrt{\frac{2}{1+\lambda}} z\right).
\end{equation*}
It follows that 
\begin{align*}
\int_{-2}^2 \frac{1}{(z-x)^{\lambda}} \left(4- x^2\right)^{\lambda-1/2}dx & \propto \nonumber = [G(z)]^{\lambda} = \left[\int_{-2}^2 \frac{1}{z-x} \sqrt{4-x^2}\frac{dx}{2\pi}\right]^{\lambda}.
 \end{align*}
Using the fact that the Wigner distribution is the Markov transform of the arcsine distribution (see \cite{Ker} p. 64), one finally gets
\begin{pro}For $\lambda > 0$,
\begin{align}\label{I1}
\int_{-2}^2 \frac{1}{(z-x)^{\lambda}} \left(4- x^2\right)^{\lambda-1/2}dx & \propto  \left[\int_{-2}^2 \frac{1}{z-x} \sqrt{4-x^2}\frac{dx}{2\pi}\right]^{\lambda}
\\& \nonumber = \exp -\lambda \int_{-2}^2 \log (z-x)\frac{1}{\pi\sqrt{4-x^2}}dx.
 \end{align}
\end{pro} 
\begin{nota}[Gauss hypergometric function]
Recall that the Gauss hypergeometric function of a complex variable $z$ is defined by 
\begin{equation*}
{}_2F_1(a,b,c;z) = \sum_{n \geq 0}\frac{(a)_n(b)_n}{(c)_n} \frac{z^n}{n!},\quad |z| < 1
\end{equation*}
where $c \in \mathbb{R} \setminus \mathbb{Z}_-$ and $(a)_n$ is the Pochhammer symbol defined by: 
\begin{equation*}
(a)_0 = 1, \quad (a)_n = (a+n-1)\cdots (a+1)a,\, n \geq 1.
\end{equation*}
When $c > b > 0$, ${}_2F_1$ admits the integral representation: 
\begin{equation*}
{}_2F_1(a,b,c;z) = \frac{\Gamma(c)}{\Gamma(c-b)\Gamma(b)} \int_0^1 (1-tz)^{-a}t^{b-1}(1-t)^{c-b-1}dt.
\end{equation*}
Then, it is an easy exercice to see that 
\begin{equation*}
\int_{-2}^2 \frac{1}{(z-x)^{\lambda}} \left(4- x^2\right)^{\lambda-1/2} \propto \frac{1}{(z+2)^{\lambda}}{}_2F_1\left(\lambda, \lambda +\frac{1}{2}, 2\lambda +1; \frac{4}{z+2}\right)
\end{equation*}
for enough large $z$. To recover (\ref{I1}), one uses the identity 
\begin{equation*}
{}_2F_1\left(\lambda, \lambda +\frac{1}{2}, 2\lambda +1; z \right) \propto \frac{1}{(\sqrt{1-z} +1)^{2\lambda}}
\end{equation*}
for $z$ in the unit disc. 
\end{nota}

\subsection{Second family}
The density of $\mu_{\lambda}$ reads
\begin{equation*}
\mu_{\lambda}(dx) \propto \left(1- \frac{x^2}{2\lambda}\right)^{\lambda-3/2} {\bf 1}_{[-\sqrt{2\lambda},\sqrt{2\lambda}]}(x)dx, \quad \lambda > 1/2,
\end{equation*} 
and we map it using $x \mapsto \sqrt{\lambda/2}x$ to 
\begin{equation*}
\left(4- x^2\right)^{\lambda-3/2} {\bf 1}_{[-2,2]}(x)dx, \quad \lambda > 1/2. 
\end{equation*}
Now, 
\begin{equation*}
f_{\lambda}^{-1}(z) = \sqrt{\frac{2}{\lambda}} G \left(\sqrt{\frac{2}{\lambda}} z\right)
\end{equation*}
so that one gets 
\begin{equation*}
\int_{-2}^2 \frac{1}{(z-x)^{\lambda}} \left(4- x^2\right)^{\lambda-3/2} dx \propto \frac{G^{\lambda}(z)}{1- G^2(z)} = \tilde{u}_{\lambda}(G)(z).  
\end{equation*}
Now recall that, for $z \in \mathbb{C} \setminus [-2,2]$,  
\begin{equation*}
G(z) = \int_{-2}^2\frac{1}{z-x} \sqrt{4-x^2}\frac{dx}{2\pi} = \frac{z - \sqrt{z^2-4}}{2} = \frac{2}{z + \sqrt{z^2-4}}, 
\end{equation*}
and that $G^2(z) + 1 = zG(z)$. Then 
\begin{equation*}
1- G^2(z) = 2 - zG(z) = 2 - \frac{2z}{z + \sqrt{z^2-4}} = \sqrt{z^2-4}G(z).
\end{equation*} 
Using (\cite{Spei}) 
\begin{equation*}
\frac{1}{\pi}\int_{-2}^2\frac{1}{z-x} \frac{dx}{\sqrt{4-x^2}} = \frac{1}{\sqrt{z^2-4}} = \exp -\frac{1}{2}\log [(z-2)(z+2)]
\end{equation*} 
for  suitable $z$, one gets 
\begin{pro}
\begin{equation}\label{I2}
\int_{-2}^2 \frac{1}{(z-x)^{\lambda}} \left(4- x^2\right)^{\lambda-3/2} dx \propto \frac{G^{\lambda-1}(z)}{\sqrt{z^2-4}} = G^{\lambda-1}(z)G_{\arcsin}(z).  
\end{equation}
Moreover,
\begin{equation*}
\tau_{\lambda}(dx) = \left(1-\frac{1}{\lambda}\right)\frac{1}{\pi} \frac{1}{\sqrt{4-x^2}}{\bf 1}_{[-2,2]}(dx) + \frac{1}{\lambda} \frac{\delta_{-2} + \delta_2}{2}(dx),  
\end{equation*}
and is not a probability measure unless $\lambda \geq 1$. 
\end{pro}
Next, for $\lambda \geq 1$, 
\begin{equation*}
T_{\lambda}(z) := G^{1-1/\lambda}(z)G_{\arcsin}^{1/\lambda}(z), \quad z \in \mathbb{C} \setminus [-2,2],
\end{equation*}
is the CST of a probability measure $\nu_{\lambda}$. This is readily checked using Lemma II. 2.2 in \cite{Sho}. More precisely, one has for $\Im(z) > 0$
\begin{equation*}
\arg[G^{1-1/\lambda}(z)G_{\arcsin}^{1/\lambda}(z)] = \left(1-\frac{1}{\lambda}\right)\arg[G(z)] + \frac{1}{\lambda}\arg[G_{\arcsin}(z)] \in ]-\pi,0[
\end{equation*}
so that $T_{\lambda}$ is of imaginary type (maps the upper half-plane into the lower half-plane), and 
\begin{equation*}
\lim_{y \rightarrow \infty} iyT_{\lambda}(iy) = \lim_{y \rightarrow \infty}  [iyG(iy)]^{1-1/\lambda}[iyG_{\arcsin}(iy)]^{1/\lambda}.
\end{equation*}
For $1/2 < \lambda < 1$, one easily gets the inequality
\begin{equation*}
2\arg[G_{\arcsin}(z)] < \arg[T_{\lambda}(z)] <  \arg[G_{\arcsin}(z)] - \arg[G(z)]
\end{equation*}
and there is no guarantee to that $\nu_{\lambda}$ is a probability distribution. 
In order to characterize $\nu_{\lambda}, \lambda \geq 1$, the Markov transform of $\tau_{\lambda}$, we will use results by Cifarelly and Regazzini (\cite{Ker} p. 51).
In fact, since $\tau_{\lambda}$ is a compactly-supported probability measure, then $\nu_{\lambda}$ is absolutely continuous with density proportional to 
\begin{equation}\label{F1}
 \sin (\pi F_{\lambda}(x)) \exp -\int_{-2}^2 \log|x - u|\tau_{\lambda}(du),
\end{equation} 
where
\begin{align*}
\pi F_{\lambda}(x) &:= \tau_{\lambda}(]-\infty,x]) 
\\&= \left\{
\begin{array}{lll}
0 & \textrm{if} & x < -2, \\
\pi/(2\lambda) & \textrm{if} & x = -2,\\
(1-1/\lambda)[\arcsin (x/2)+ \pi/2] & \textrm{if} & x \in [-2,2[,\\
\pi & \textrm{if} & x \geq 2 
\end{array}
\right.. \end{align*}
Note that since (\ref{F1}) is valid when $\tau_{\lambda}$ is the arcsine distribution and $\nu_{\lambda}$ is  the Wigner distribution, one deduces that 
\begin{equation*}
\exp -\int_{-2}^2 \log|x - u|\frac{1}{\pi} \frac{1}{\sqrt{4-u^2}}du, \quad x \in [-2,2],
\end{equation*}
does not depend on $x$ (is constant). This striking result may be used to derive the density of $\nu_{\lambda}$ in our case and in the forthcoming ones since $\tau_{\lambda}$ is a convex linear combination of the arcsine distribution and a discrete probability measure. In the case in hand, easy computations yield 
\begin{pro}
\begin{equation*}
\frac{\nu_{\lambda}(dx)}{dx} \propto \cos \left[\left(1-\frac{1}{\lambda}\right)\arcsin \frac{x}{2}\right] \frac{1}{(4-x^2)^{1/2\lambda}}{\bf 1}_{]-2,2[}(x), \lambda \geq 1.
\end{equation*}
\end{pro}
Note that $\lambda=1$ corresponds to the arcsine distribution while $\lambda = \infty$ correponds to the Wigner distribution. Thus, $\nu_{\lambda}$ interpolates between them. 
The reader may wonder how to compute the normalizing constant or the moments of $\nu_{\lambda}$. This will be clear after dealing with the two remaining families $\mu_{\lambda}$.

\begin{nota}
One can derive (\ref{I2}) using the Gauss hypergeometric function. This time, 
\begin{equation*}
\int_{-2}^2 \frac{1}{(z-x)^{\lambda}} \left(4- x^2\right)^{\lambda-3/2} \propto \frac{1}{(z+2)^{\lambda}}{}_2F_1\left(\lambda, \lambda - \frac{1}{2}, 2\lambda -1; \frac{4}{z+2}\right)
\end{equation*}
for enough large $z$, and one makes use of 
\begin{equation*}
{}_2F_1\left(\lambda - \frac{1}{2}, \lambda, 2\lambda-1; z \right) \propto \frac{1}{\sqrt{1-z}(\sqrt{1-z} +1)^{2(\lambda-1)}}
\end{equation*}
for $z$ in the unit disc since
\end{nota}

\section{Markov transforms: non symmetric measures}
\subsection{Third family}
The probability distribution $\mu_{\lambda}$ has the density 
\begin{equation*}
\left(1 - \frac{\sqrt{2\lambda-1}x -1}{2\lambda}\right)^{\lambda-1/2}\left(1 + \frac{\sqrt{2\lambda-1}x -1}{2\lambda}\right)^{\lambda-3/2}
\end{equation*}
where
\begin{equation*}
x \in  \left[\frac{1-2\lambda}{\sqrt{2\lambda-1}}, \frac{1+2\lambda}{\sqrt{2\lambda-1}}\right], \quad \lambda > 1/2.
\end{equation*}
Moreover
\begin{equation*}
f_{\lambda}(z) = \frac{\lambda^2}{2\lambda-1}z + \frac{1}{\sqrt{2\lambda-1}} + \frac{1}{z} = K_{1/\sqrt{2\lambda-1}, \lambda/ \sqrt{2\lambda-1}}(z). 
\end{equation*}
Thus,
\begin{align*}
f_{\lambda}^{-1}(z) = G_{1/\sqrt{2\lambda-1}, \lambda/ \sqrt{2\lambda-1}}(z) &= \frac{\sqrt{2\lambda-1}}{\lambda} 
G\left[\frac{\sqrt{2\lambda-1}}{\lambda}\left(z-\frac{1}{\sqrt{2\lambda-1}}\right)\right]
\\& =  \frac{\sqrt{2\lambda-1}}{\lambda} G\left[\frac{\sqrt{2\lambda-1}z-1}{\lambda}\right].
\end{align*}
Now, the image of $\mu_{\lambda}$ under the map $x \mapsto (\sqrt{2\lambda-1}x-1)/\lambda$ transforms its density to 
\begin{equation*}
 \left(1 -\frac{x}{2}\right)^{\lambda-1/2}\left(1 + \frac{x}{2}\right)^{\lambda-3/2}{\bf 1}_{[-2,2]}(x)
 \end{equation*}
and one easily sees that 
\begin{equation*}
\int_{-2}^2\frac{1}{(z-x)^{\lambda}} \left(4-x^2\right)^{\lambda-3/2}\left(2- x\right)dx \propto \frac{G^{\lambda}(z)}{1- G(z)} := \tilde{u}_{\lambda}(G)(z).
 \end{equation*} 
Now note that $(1-G(z))^2 = 1+G^2(z) -2G(z) = (z-2)G(z)$ which yields $1- G(z) = \sqrt{z-2}\sqrt{G(z)}$ where the branch of the square root is taken so that $1-G$ is of positive imaginary type (i.e. maps the upper half plane into itself). As a result, 
\begin{pro}
\begin{equation*}
\int_{-2}^2\frac{1}{(z-x)^{\lambda}} \left(4-x^2\right)^{\lambda-3/2}\left(2- x\right)dx \propto \frac{G^{\lambda}(z)}{1- G(z)} = \frac{G^{\lambda-1/2}(z)}{\sqrt{z-2}}.
\end{equation*}
\end{pro}
In this case $\tau_{\lambda}$ satisfies for suitable $z$ 
\begin{equation*}
\int_{\mathbb{R}} \log (z-x) \tau_{\lambda}(dx) = -\left(1-\frac{1}{2\lambda}\right) \log(G(z)) - \frac{1}{2\lambda} \log\frac{1}{z-2} 
\end{equation*} 
whence we deduce that  
\begin{equation*}
\tau_{\lambda}(dx) = \left(1-\frac{1}{2\lambda}\right)\frac{1}{\pi} \frac{1}{\sqrt{4-x^2}}{\bf 1}_{[-2,2]}(dx) + \frac{1}{2\lambda}\delta_2(dx) 
\end{equation*}
which is a probability measure for all $\lambda > 1/2$. Besides, the same arguments used before show that for suitable $z$
\begin{equation*}
z \mapsto \frac{G^{1-1/(2\lambda)}(z)}{(z-2)^{1/2\lambda}},
\end{equation*}
is the CST of a probability distribution $\nu_{\lambda}, \lambda > 1/2$ which is absolutely continuous of density given by
\begin{pro}
\begin{equation*}
\frac{\nu_{\lambda}(dx)}{dx} \propto \sin \left[\left(1-\frac{1}{2\lambda}\right)\left(\arcsin \frac{x}{2} + \frac{\pi}{2}\right)\right] \frac{1}{(2-x)^{1/2\lambda}}{\bf 1}_{]-2,2[}(x), \lambda > 1/2.
\end{equation*}
\end{pro}

\subsection{Fourth family}
The density of $\mu_{\lambda}$ is given by
\begin{equation*}
\left(1 - \frac{\sqrt{2\lambda-1}x +1}{2\lambda}\right)^{\lambda-3/2}\left(1 + \frac{\sqrt{2\lambda-1}x +1}{2\lambda}\right)^{\lambda-1/2}
\end{equation*}
where
\begin{equation*}
x \in  \left[\frac{-1-2\lambda}{\sqrt{2\lambda-1}}, \frac{-1+2\lambda}{\sqrt{2\lambda-1}}\right], \quad \lambda > 1/2.
\end{equation*}
The density of the image of $\mu_{\lambda}$ under the map $x \mapsto (\sqrt{2\lambda-1}x+1)/\lambda$ reads 
\begin{equation*}
 \left(1 -\frac{x}{2}\right)^{\lambda-3/2}\left(1 + \frac{x}{2}\right)^{\lambda-1/2}{\bf 1}_{[-2,2]}(x) \propto \left(4-x^2\right)^{\lambda-3/2}\left(2 +x\right){\bf 1}_{[-2,2]}(x). 
 \end{equation*}

The same scheme used to deal with the third family gives 
\begin{equation*}
\tau_{\lambda}(dx) = \left(1-\frac{1}{2\lambda}\right)\frac{1}{\pi} \frac{1}{\sqrt{4-x^2}}{\bf 1}_{[-2,2]}(dx) + \frac{1}{2\lambda}\delta_{-2}(dx) 
\end{equation*}
and that the CST of $\nu_{\lambda}$ is given by 
\begin{equation*}
\frac{G^{1-1/(2\lambda)}(z)}{(z+2)^{1/2\lambda}}.
\end{equation*}
However, the density of $\nu_{\lambda}$ is somewhat different from the one in the previous case: 
\begin{equation*}
\frac{\nu_{\lambda}(dx)}{dx} \propto \cos \left[\left(1-\frac{1}{2\lambda}\right)\arcsin \frac{x}{2}\right] \frac{1}{(2+x)^{1/2\lambda}}{\bf 1}_{]-2,2[}(x), 
\lambda > 1/2.
\end{equation*}

\begin{note}
1/For the four families, the measure $\tau_{\lambda}$ may be mapped to 
\begin{equation*}
 \frac{1}{\pi} \frac{1}{\sqrt{1-x^2}}{\bf 1}_{[-1,1]}(dx) + M_{\lambda}\delta_{-1}(dx) + N_{\lambda} \delta_{1}(dx)
 \end{equation*} 
for some positive constants $M_{\lambda},N_{\lambda}$. A more wider class of measures including the above one were considered in \cite{Koor}. \\ 
2/For $0< a < 1, b =1-a$, one can always define an operation 
\begin{equation*}
(\nu_1,\nu_2) \mapsto \nu / G_{\nu} = G_{\nu_1}^aG_{\nu}^b,
\end{equation*}
where $\nu_1,\nu_2,\nu$ are probability measures. However, every probability measure will be idompotent and the operation is not commutative unless $a=b=1/2$.    
\end{note}

\section{On the moments of the second Markov transform}
It is known that the Wigner distribution is a universal limiting distribution: it is the spectral distribution of large rescaled random matrices from the so called Wigner ensemble, the limiting distribution of the rescaled Plancherel transition of the growth process for Young diagrams (\cite{Ker1}) and more generally of the rescaled transition measure of rectangular diagrams associated with roots of some adjacent orthogonal polynomials (\cite{Ker2}). It also plays a crucial role in free probability theory where it appears as the central limiting distribution of the sum of free random variables (\cite{Spei}, \cite{Bia}). Note also that the standard arcsine distribution is the central limiting distribution for the $t$-convolution with $t=1/2$ 
(\cite{Bo}, \cite{BKW}, \cite{Kry}) and is the limiting distribution of the so-called Shrinkage process (\cite{Ker1}).\\
As a matter a fact, it is interesting to find a parallel to the above facts when the Wigner law is replaced by the Markov transform $\nu_{\lambda}$ corresponding to the second family. More interesting is to define a convolution operation that interpolates the $t$-convolutions for $t=1/2$ ($\lambda=1$) and $t=1$ ($\lambda = \infty$, free convolution) having $\nu_{\lambda}$ as a central limiting distribution (\cite{Cab}). Since $\nu_{\lambda}$ is symmetric and compactly-supported, it is entirely determined by its even moments and we claim

\begin{pro}
The  normalizing constant of $\nu_{\lambda}$ is given by 
\begin{equation*}
c_{\lambda} = 2^{1-1/\lambda}\sqrt{\pi}\frac{\Gamma(1-1/(2\lambda))\Gamma(3/2-1/(2\lambda))}{\Gamma(2-1/\lambda)}
\end{equation*}
and the even moments of $\nu_{\lambda}$ may be expressed as
 \begin{equation*}
m_{2n}^{\lambda}:= \int_{-2}^2 x^{2n}\nu_{\lambda}(dx)  = 2^{2n} {}_3F_2\left(\substack{\displaystyle -n,1-1/(2\lambda),3/2 - 1/(2\lambda)\\ \displaystyle 2-1/\lambda, 1}; 1\right).
\end{equation*}
Moreover, the moments are polynomials in the variable $(1/\lambda)$.
\end{pro}
{\it Proof}: make the change de variable $x \mapsto 2\sin x$ in the integral
\begin{equation*}
 \int_{-2}^{2} x^{2n} \cos \left[\left(1-\frac{1}{\lambda}\right)\arcsin\frac{x}{2}\right] [4-x^2]^{-1/(2\lambda)} dx, \, n \geq 0, 
 \end{equation*}
to obtain
\begin{equation*}
 2^{2n+2-1/\lambda} \int_0^{\pi/2} [\sin x]^{2n} \cos \left[\left(1-\frac{1}{\lambda}\right)x\right] [\cos x]^{1-1/\lambda} dx, \, n \geq 0. 
\end{equation*}
Then, expand  
\begin{equation*}
[\sin x]^{2n} = (1-\cos^2x)^n = \sum_{k=0}^n \binom{n}{k} (-1)^k [\cos x]^{2k}
\end{equation*}
and use the formula (see \cite{Sum}, p. 177) 
\begin{equation}\label{Cauchy}
\int_0^{\pi/2}\cos [(p-q)x] [\cos x]^{p+q-2} dx = \frac{\pi}{2^{p+q-1}}\frac{\Gamma(p+q-1)}{\Gamma(p)\Gamma(q)},\, p+ q > 1,
\end{equation}
with $p+q = 2k + 3 - 1/\lambda, p-q = 1-1/\lambda$, to get 
\begin{equation*}
  \int_{-2}^{2} x^{2n} \cos \left[\left(1-\frac{1}{\lambda}\right)\arcsin\frac{x}{2}\right] [4-x^2]^{-1/2\lambda} dx =  2^{2n} \pi\sum_{k=0}^n \binom{n}{k} \frac{\Gamma(2k+2 -1/\lambda)}{\Gamma(k+2-1/\lambda)k!} \left(-\frac{1}{4}\right)^k.
\end{equation*}
This may be expressed through ${}_3F_2$ hypergeometric series as follows: write the binomial coefficient as 
\begin{equation*}
\binom{n}{k} = \frac{n!}{k!(n-k)!} = (-1)^k\frac{(-n)_k}{k!}
\end{equation*}
and use the duplication formula to rewrite
\begin{equation*}
\sqrt{\pi} \Gamma\left(2k + 2 -\frac{1}{\lambda}\right) = 2^{2k+1-1/\lambda}\Gamma\left(k+1 -\frac{1}{2\lambda}\right)\Gamma\left(k+3/2 - \frac{1}{2\lambda}\right).
\end{equation*} 
It follows that 
\begin{align*}
 \int_{-2}^{2} x^{2n} \cos \left[\left(1-\frac{1}{\lambda}\right)\arcsin\frac{x}{2}\right] [4-x^2]^{-1/2\lambda} dx = c_{\lambda}2^{2n}
 {}_3F_2\left(\substack{\displaystyle -n,1-1/(2\lambda),3/2 - 1/(2\lambda)\\ \displaystyle 2-1/\lambda, 1}; 1\right).
\end{align*}
To prove the last claim, let $y = 1/\lambda$ and expand:
\begin{align*}
\left(1-\frac{y}{2}\right)_k &= \left(k-\frac{y}{2}\right)\left(k-1-\frac{y}{2}\right) \cdots \left(1-\frac{y}{2}\right)
\\& = \left(-\frac{1}{2}\right)^k (y-2)(y-4)\cdots(y-2k+2)(y-2k),
\end{align*}
and similarly 
\begin{eqnarray*}
\left(\frac{3}{2}-\frac{y}{2}\right)_k &=& \left(-\frac{1}{2}\right)^k (y-3)(y-5)\cdots(y-2k-1),\\
\left(2-y\right)_k &=& (-1)^k(y-2)(y-3)\cdots (y-k-1).
\end{eqnarray*}
The proof ends after forming the ratio
\begin{equation*}
\frac{(1-y/2)_k(3/2 - y/2)_k}{(2-y)_k} = \left(-\frac{1}{4}\right)^k (y-k-2)\cdots (y- 2k)(y-2k-1) 
\end{equation*}
so that 
\begin{equation*}
m_{2n}^{\lambda} := p_n(y) = \sum_{k=0}^n \binom{n}{k}4^{n-k} \frac{(y-k-2)\cdots (y- 2k)(y-2k-1)}{k!}. \qquad \blacksquare
\end{equation*} 
\begin{nota}
The exponential generating series of $(p_n)_n$ is given by 
\begin{align*}
\sum_{n \geq 0}p_n(y) \frac{z^n}{n!} &= \sum_{0 \leq k \leq n} \frac{4^{n-k}}{(n-k)!k!} \frac{(y-k-2)\cdots (y- 2k)(y-2k-1)}{k!} z^n
\\& = e^{4z}\sum_{k\geq 0} \frac{(y-k-2)\cdots (y- 2k)(y-2k-1)}{k!} \frac{z^k}{k!} 
\\& = 2^{1-y}e^{4z}{}_2F_2\left(\substack{\displaystyle 1-y/2, (3-y)/2,\\ \displaystyle 2-y, 1} ; 4z\right).  
\end{align*}

\end{nota}
\subsection{Combinatorics}
For $\lambda = 1$, one recovers the moments of the arcsine distribution given by (use duplication formula)
\begin{equation*}
\binom{2n}{n} = 2^{2n} \frac{(1/2)_n}{(1)_n} = 2^{2n}{}_2F_1\left(\substack{\displaystyle -n, 1/2\\ \displaystyle 1}; 1\right).
\end{equation*}
When $\lambda = \infty$, one recovers the moments of the Wigner distribution: 
\begin{equation*}
\frac{1}{n+1}\binom{2n}{n} = 2^{2n} \frac{(1/2)_n}{(2)_n} = 2^{2n}{}_2F_1\left(\substack{\displaystyle -n, 3/2\\ \displaystyle 2}; 1\right).
\end{equation*}
Both moments are integers and it is known that they count the shifted and the ordinary Catalan paths respectively (see \cite{Ker}, p. 64). Unfortunately, $m_{2n}^{\lambda}$ is not positive integer-valued in general therefore no similar combinatorial interpretation can be given for all values of $\lambda$. For instance, when $\lambda = 2$, one gets 
\begin{equation*}
m_{2n}^2 = 2^{2n+1} \sum_{k=0}^n (-1)^k \binom{n}{k} \frac{1}{\sqrt{\pi}} \frac{\Gamma(2k+3/2)}{\Gamma(2k+2)} 
= 2^{2n-1} \sum_{k=0}^n (-1)^k \binom{n}{k}\binom{4k+2}{2k+1} \frac{1}{2^{4k}}  
\end{equation*}
by Gauss duplication formula so that the few first even moments are $3/2, 31/8, 187/16, 4859/128$.  

\begin{nota}[Non symmetric Markov transforms]
The computations of the moments of the non symmetric Markov transforms need more care. Consider for instance the one corresponding to the last family, then one has to compute
after the change of variables $x \mapsto 2\sin x$
\begin{equation*}
2^{n+1-1/(2\lambda)}\int_{-\pi/2}^{\pi/2}[\sin x]^{n}\cos\left[\left(1-\frac{1}{2\lambda}\right)x\right] \frac{\cos x}{(1+\sin x)^{1/2\lambda}}dx.
\end{equation*}
To proceed, use the binomial Theorem to expand
\begin{equation*}
\frac{1}{(1+\sin x)^{1/2\lambda}} = \sum_{k \geq 0}\frac{(1/2\lambda)_k}{k!} (-1)^k [\sin x]^k 
\end{equation*}
and Fubini's Theorem to exchange both integral and sum signs. In fact,  
\begin{align*}
\int_{-\pi/2}^{\pi/2} |\sin x|^{n}&\cos\left[\left(1-\frac{1}{2\lambda}\right)x\right] \cos x \sum_{k \geq 0}\frac{(1/2\lambda)_k}{k!}|\sin x|^kdx 
\leq \int_{-\pi/2}^{\pi/2} \frac{\cos x}{(1-|\sin x|)^{1/2\lambda}}dx 
\\& = 2  \int_0^{\pi/2} \frac{\cos x}{(1-\sin x)^{1/2\lambda}}dx  = -\frac{[(1-\sin x)^{1-1/2\lambda}]_0^{\pi/2} }{1-1/(2\lambda)}< \infty
\end{align*}
for $\lambda > 1/2$. Now, it only remains to compute the integral 
\begin{equation*}
\int_{-\pi/2}^{\pi/2}[\sin x]^{n+k}\cos\left[\left(1-\frac{1}{2\lambda}\right)x\right]\cos xdx
\end{equation*}
which is zero when $n+k$ is odd and is evaluated similarly as the even moments of the above symmetric Markov transform when $n+k$ is even using (\ref{Cauchy}).\\
For the third family, make the substitution $x \mapsto \arcsin(x/2) + \pi/2$ in the integral
\begin{equation*}
\int_{-2}^2 x^n \sin \left[\left(1-\frac{1}{2\lambda}\right)\left(\arcsin \frac{x}{2} + \frac{\pi}{2}\right)\right] \frac{1}{(2-x)^{1/2\lambda}}dx
\end{equation*}
to get 
\begin{equation*}
2^{n+1-1/(2\lambda)} (-1)^n\int_{0}^{\pi} [\cos x]^n \sin \left[\left(1-\frac{1}{2\lambda}\right)x \right] \frac{\sin x}{(1+\cos x)^{1/2\lambda}}dx.
\end{equation*}
Then use the trigonometric formula 
\begin{equation*}
\sin a \sin b = \frac{1}{2}[\cos(a-b) - \cos(a+b)]
\end{equation*}
to split the last integral to 
\begin{equation*}
\frac{1}{2}\int_{0}^{\pi} [\cos x]^{n+1} \frac{1}{(1+\cos x)^{1/2\lambda}}dx - \frac{1}{2}\int_{0}^{\pi} [\cos x]^n \cos \left[\left(2-\frac{1}{2\lambda}\right)x \right] \frac{1}{(1+\cos x)^{1/2\lambda}}dx.
\end{equation*}
Finally, use the binomial Theorem and (\ref{Cauchy}). 
\end{nota}

\subsection{Inverses of CST for $\lambda = 2,3,4$}
Let $\lambda =n \geq 1$ be a positive integer. Recall that 
\begin{equation*}
[G_{n}(z)]^n := \left[\int_{\mathbb{R}}\frac{1}{(z-x)}\nu_{n}(dx)\right]^{n} = \tilde{u}_{n}(G(z)),
\end{equation*} 
where as before $G$ is the CST of the Wigner distribution. For the first family, $G_n = G$ for all $n$ so that $K_n = K$ and $K(z) = z+1/z$ for $z$ in a neighborhood of zero. For the remaining families, the inverse of $G_n$ is given by 
\begin{equation*}
G_n^{-1}(z) := K_n(z) = K(\tilde{u}_n^{-1}(z^n)), \qquad K := G^{-1},
\end{equation*}    
for small $|z|$, subject to the condition
\begin{equation*}
K_n(z) = \frac{1}{z} + R_n(z) =  \frac{1}{z} + \sum_{k \geq 0}r_k z^k
\end{equation*}   
where $R_n$ is an entire function known in free probability theory as the \emph{free cumulants generating function} of $\nu_n$, $(r_k)_k$ is the sequence of the free cumulants, $r_0, r_1$ are the mean and the variance of $\nu_n$ respectively.  For the second family, one has to invert 
\begin{equation*}
z \mapsto \tilde{u}_n(z) = \frac{z^n}{1-z^2}
\end{equation*} 
for small $|z|$ in the lower half unit disc (image of $G$). We are thus led to find one root of the polynomial $z^n +z^2w - w$ for complex numbers $z,w$ in a neighborhood of zero. This task is very complicated and even more, Galois theory tells us that the roots cannot be expressed by means of radicals unless $n\leq 4$. For $n=2$, easy computations show that 
\begin{equation*}
K_2(z) = \frac{2z^2+1}{z\sqrt{1+z^2}} = \frac{1}{z} + \sum_{k \geq 0}\frac{(1/2)_k}{(k+1)!}(-1)^k (k+3/2)z^{2k+1}.
\end{equation*}
Thus one sees that $r_0 = 0$, $r_1 = 3/2$ which agrees with our above computations, and that 
\begin{equation*}
r_{2k+1} = (-1)^k\left[\frac{(1/2)_k}{(1)_k}  + \frac{1}{2}\frac{(1/2)_k}{(2)_k}\right].
\end{equation*}
Thus, 
\begin{equation*}
(-1)^k2^{2k+1}r_{2k+1} = 2^{2k+1}\frac{(1/2)_k}{(1)_k} + 2^{2k}\frac{(1/2)_k}{(2)_k} = 2 \binom{2k}{k} + \frac{1}{k+1}\binom{2k}{k}.
\end{equation*}

For $n=4$, one has a quadratic polynomial and setting $v = z^2$, one is led to $v^2 + wv - w = 0$ so that
\begin{equation*}
v = \frac{-w+ \sqrt{w^2+4w}}{2}.
\end{equation*}
The case $z=3$ is more complicated and we supply one way to adress it: make the substitution $z = Z - w/3$ for suitable $Z$ to get 
\begin{equation*}
z^3 + wz - w = Z^3 - \frac{w^2}{3}Z + \frac{2}{27}w^3 - w.
\end{equation*} 
The last polynomial has the same form as  
\begin{equation*}
(a+b)^3 -  3ab(a+b) - (a^3+b^3) = 0 
\end{equation*}
which hints to look for a root of the form $Z = a+b$ where 
\begin{equation*}
ab = \frac{w^2}{3},\quad a^3+b^3 = w\left[\frac{2}{27}w^2 - 1\right]
\end{equation*}
and computations are left to the curious reader.
\begin{nota}
The above line of thinking remains valid for the non symmetric Markov transforms for which
\begin{equation*}
 \tilde{u}_n(z) = \frac{z^n}{1\pm z}.
 \end{equation*}
For $\lambda = 3$, one already has the appropriate form of the polynomial and there is no need to make the above substitution. However, the case $\lambda = 4$ needs more developed techniques since the polynomial is of degree four and is not quadratic.  
\end{nota}

{\bf Acknowledgments}: this work was fully supported by CRC 701. The author finished the paper while visiting Georgia Institute of Technology and University of Virginia. He actually wants to thank H. Matzinger, C. Houdre and the administrative staff for their hospitality and their help with references. A special thank is given to Professor C. F. Dunkl for fruitful discussions held at UVa and for his help while doing tedious computations, in particular for numerical simulations.


\begin{thebibliography}{99}
\bibitem{Bia}\emph{P. Biane}. Representations of symmetric groups and free probability. {\it Adv. Math.} {\bf 138}. no. 1. 1998, 126-181.
\bibitem{Ber}\emph{H. Bercovici, D. V. Voiculescu}. L\'evy-Hincin-type theorems for multiplicative and additive free convolutions. {\it Pacific. J. Math}. {\bf 153}. 1992, 217-248.
\bibitem{Bo}\emph{M. Bozejko}. Deformed Free Probability of Voiculescu, {\it R. I. M. S. Kokyuroku}. {\bf 1227}. 2001, 96Ð113. 
\bibitem{Boz}\emph{M. Bozejko, N. Demni}. Generating functions of Cauchy-Stieltjes type for orthogonal polynomials. {\it To appear in Infinite. Dimen. Anal. Quantum Probab. Relat. Top.}
\bibitem{BKW}\emph{M. Bozejko, A. D. Krystek, L. J. Wojakowski}. Remarks on the $r$- and .-convolutions, {\it M. Zeit}. {\bf 253}. 2006, 177Ð196. 
\bibitem{Dem}\emph{N. Demni}. Ultraspherical type generating functions for orthogonal polynomials. {\it Submitted to Probab. Math. Statist. Available on arXiv}.
\bibitem{Cab}\emph{T. Cabanal-Duvillard}. Un th\'or\`eme central limite pour des variables al\'atoires non-commutatives. {\it C. R. A. S. t. 325}, {\bf SŽrie I}. 1997, 1117-1120. 
\bibitem{Hirs}\emph{I. I. Hirschmann, D. V. Widder}. Generalized inversion formulas for convolution transforms. II. {\it Duke Math. J.} {\bf 17}. 391-402. 
\bibitem{Ker}\emph{S. Kerov}. Interlacing measures. {\it Amer. Math. Soc. Transl. Ser 2}. {\bf 181}.  1998, 35-83. 
\bibitem{Ker1}\emph{S. Kerov}.Transition probabilities of continuous Young diagrams and Markov's moment problems. {\it Func. Anal. Appl}. {\bf 27}. no. 2. 1993, 32-49. 
\bibitem{Ker2}\emph{S. Kerov}. Asymptotic separation of roots of orthogonal polynomials. {\it Algebra and Analysis}. {\bf 5}. no .5. 1993, 68-86. 
\bibitem{Koor}\emph{T. Koornwinder}. Orthogonal polynomials with weight function $(1-x)^{\alpha}(1+x)^{\beta} + M\delta(x+1) + N\delta(x-1)$. {\it Canad. Math. Bull}. {\bf 27}. no. 2. 1984, 205-214.
\bibitem{Kry}\emph{A. D. Krystek, H. Yoshida}. The combinatorics of the $r$-free convolution. {\it Infinite Dimen. Anal. Quantum. Probab}. {\bf 6}. no. 4. 2003, 619-627. 
\bibitem{Sch}\emph{J. H. Schwarz}. The generalized Stieltjes transform and its inverse. {\it J. Math. Phys.} {\bf 46}. no. 1. 2005, 013501, 8pp.
\bibitem{Sho}\emph{J. A. Shohat, J. D. Tamarkin.} The Problem of Moments. {\it Amer. Math. Soc.} NY, 1943.
\bibitem{Spei}\emph{R. Speicher}. Combinatorics of Free Probability Theory. {\it Lectures. I. H. P. Paris}. 1999.
\bibitem{Sum}\emph{D. B. Sumner}. An inversion formula for the generalized Stieltjes transform. {\it Bull. Amer. Math. Soc}. {\bf 55}. 1949, 174-183. 
\end{thebibliography}
\end{document}